\documentclass[11pt]{article}
\usepackage{graphicx}
 \usepackage{mathptmx,color}      
\usepackage[titletoc,title]{appendix}
\usepackage[numbers]{natbib}
\usepackage[T1]{fontenc}
\usepackage[applemac]{inputenc}
\usepackage[all]{xy}
\usepackage{enumerate}
\usepackage{comment}
\usepackage{amsmath,amssymb,stmaryrd,bm}
\usepackage[british]{babel}
\usepackage{a4wide}
\usepackage[sans]{dsfont}
\usepackage{amsfonts}
\usepackage{amsthm}
\usepackage[mathscr]{euscript}
\usepackage{graphicx}
\usepackage{mathrsfs}
\usepackage{latexsym}
\usepackage{setspace}
\usepackage{amsthm,amsmath}

\newcommand{\p}{^}

\newcommand{\ol}{\overline}

	\newcommand{\rmd}{\mathrm{d}}

\newcommand{\Bscr}{\mathscr{B}}

\newcommand{\Escr}{\mathscr{E}}
\newcommand{\Fscr}{\mathscr{F}}

\newcommand{\Hscr}{\mathscr{H}}

\newcommand{\Kscr}{\mathscr{K}}
\newcommand{\Lscr}{\mathscr{L}}
\newcommand{\Mscr}{\mathscr{M}}

\newcommand{\Pscr}{\mathscr{P}}

\newcommand{\zt}{\zeta}
\newcommand{\et}{\eta}

\newcommand{\ep}{\varepsilon}

\newcommand{\sig}{\sigma}

\newcommand{\om}{\omega}
\newcommand{\Om}{\Omega}

\newcommand{\Ebb}{\mathbb{E}}

\newcommand{\Fbb}{\mathbb{F}}

\newcommand{\Pbb}{\mathbb{P}}

\newcommand{\Rbb}{\mathbb{R}}

\newcommand{\tr}{\mathrm{tr}}

\newcommand{\mub}{\mu}

\newtheorem{theorem}{Theorem}[section]
\newtheorem{lemma}[theorem]{Lemma}

\theoremstyle{definition}

\newtheorem{assumption}[theorem]{Assumption}

\newtheorem{remark}[theorem]{Remark}

\numberwithin{equation}{section}

\makeatletter
\def\timenow{\@tempcnta\time
\@tempcntb\@tempcnta
\divide\@tempcntb60
\ifnum10>\@tempcntb0\fi\number\@tempcntb
:\multiply\@tempcntb60
\advance\@tempcnta-\@tempcntb
\ifnum10>\@tempcnta0\fi\number\@tempcnta}
\makeatother

\title{Degenerate McKean--Vlasov SDEs\\ with Yamada--Osgood Coefficients\\ and Arbitrary  Modulus of Continuity\\ in Wasserstein Distance}
\author{Rainer Buckdahn\footnote{Shandong University, Weihai and Universit\'e de Bretagne Occidentale, Brest. Contact:  rainer.buckdahn@univ-brest.fr}, Paolo Di Tella\footnote{Technische Universit\"at Dresden, Contact: paolo.di\_tella@tu-dresden.de} and Juan Li\footnote{Shandong University, Weihai, Contact: juanli@sdu.edu.cn}   \\  
     Current version: \today, \timenow }
\date{}

\begin{document}
\maketitle

\begin{abstract}
We establish the well-posedness of strong solutions for a class of degenerate mean-field stochastic differential equations with coefficients satisfying non-Lipschitz conditions. The diffusion coefficient is not required to be non-degenerate. Both the diffusion and the drift coefficients may depend on the law through the Kantorovich-Rubinstein metric with an arbitrary modulus of continuity. The drift and the diffusion satisfy a dimension-dependent Osgood-type condition, that becomes particularly sharp in dimension one. The proof combines the classical Yamada--Watanabe argument with Schauder's fixed-point theorem and does not rely on Gaussian estimates, allowing us to cover degenerate diffusions where such estimates are not available. This fills a gap in the existing literature.
\end{abstract}

\section{Introduction}\label{sec:introduction}

In the present paper, we study the existence of a strong solution of the \(\mathbb{R}^m\)-dimensional mean field stochastic differential equation (SDE)

\[
X_t = \xi + \int_0^t b(s, X_s, \mathcal{L}(X_s)) \, \mathrm{d}s + \int_0^t \sigma(s, X_s, \mathcal{L}(X_s)) \, \mathrm{d}B_s, \qquad t \in [0,T],
\]
for continuous coefficients of sublinear growth satisfying the conditions

\[
\|\sigma(t,x,\mu)-\sigma(t,y,\nu)\|_F \leq h(\|x-y\|) + \chi(W_1(\mu,\nu)),
\]
\[
\|b(t,x,\mu)-b(t,y,\nu)\| \leq g(\|x-y\|) + \chi(W_1(\mu,\nu)),
\]
where \(h\) and \(\chi\) are strictly increasing continuity moduli such that \(h^2\) satisfies the Osgood condition in \([0,1]\). If \(m=1\), \(g\) is assumed to be continuous, concave, strictly increasing and satisfying an Osgood condition too. If \(m>1\), these assumptions are adapted to a dimension-dependent condition that is stronger than for $m=1$(see Assumption \ref{ass:ass.coeff} for more details). We stress that we do not require any degeneracy conditions on the diffusion coefficient \(\sigma\), and that the regularity in Wasserstein distance can be very rough, in fact, \(\chi\) can be any increasing modulus of continuity with \(\chi(0+)=0\).

While recent results achieve well-posedness for equations with more irregular coefficients than those studied in the present paper, they often rely on a specific structure such as the uniform ellipticity of the diffusion matrix. For instance, we recall the very general results obtained in \cite{MV20}, where, for sublinear-growth coefficients with a linear dependence on the measure, Mishura and Veretennikov have shown the existence of a weak solution. In Issoglio et al.~\cite{IPRT26}, the authors have established the weak well-posedness of a class of degenerate measure-dependent kinetic SDEs with a distributional drift in an anisotropic negative Besov space. In \cite{PR25}, Pascucci and Rondelli have shown the existence of a unique strong solution for bounded and globally H\"older continuous coefficients satisfying an \(\alpha\)-H\"older-Wasserstein distance bound in the measure. In \cite{PRV26}, Pascucci et al.\ have generalized the results of \cite{PR25} to continuous sublinear coefficients, also with relaxing (but not completely dropping) the non-degeneracy condition on the diffusion matrix.

We stress that although all these papers allow rather general coefficients, they do not cover the case of an irregular diffusion matrix that may vanish at zero, i.e., a truly degenerate diffusion, while simultaneously allowing rough measure dependence in the weak Wasserstein metric and non-Lipschitz spatial coefficients. The main contribution of this paper is to fill this gap. 

Our approach combines the classical Yamada--Watanabe argument (see \cite{YW711, YW712, AS13}) with the Schauder fixed-point theorem. The main idea is to use this method here, not to show the strong uniqueness of the solution--it is proven separately, under different assumptions-- as classically done in the literature, but rather to show the continuity of the mapping $\Phi(\mu)_t=\Lscr(X_t)$: This is the key step to apply Shouder's theorem. Unlike \cite{PR25}, as we do not rely on Gaussian estimates for the fundamental solution of the associated PDE, we can handle the case of degenerate diffusions where such estimates are not available. 

We also point out that the study of mean-field SDEs without ellipticity conditions has already been considered in the literature. For example, in \cite{LM16}, Li and Min studied weak well-posedness of certain path-dependent mean-field SDEs under very weak assumptions on the drift, assumed therein to be only bounded and measurable in space and to have a continuity modulus in measure. However, the diffusion matrix in \cite{LM16} is assumed to be invertible with bounded inverse and, more importantly, independent of the law. In contrast, our diffusion coefficient may be degenerate and may depend on the measure in a rough way.

In the present paper we also show under additional conditions (see Theorem \ref{thm:un}) the uniqueness of the strong solution. For this end we follow \cite{LeeShih04} and prove strong uniqueness by weakening the classical monotonicity condition, adapting it to our non-Lipschitz framework.

\section{Existence of a Strong Solution}\label{sec:existence}
Let an $\Rbb\p d$-valued standard Brownian motion $B$ be given on a complete probability space $(\Omega,\Fscr,\Pbb)$ and denote by $\Fbb\p B$ the completed Brownian filtration. Let $\sig: [0,T]\times\Rbb\p m\times\Pscr_1(\Rbb\p m)\longrightarrow\Rbb\p{m\times d}$ and  $b: [0,T]\times\Rbb\p m\times\Pscr_1(\Rbb\p m)\longrightarrow\Rbb\p{m}$ be jointly continuous functions. Here, $\Pscr_p(\Rbb\p m)$, $p\geq1$, denotes the space of probability measures $\mu$ on $\Rbb\p m$ endowed with the  Borel-$\sigma$-algebra $\Bscr(\Rbb\p m)$ such that $\int_{\Rbb\p m}\|x\|\p p\mu(\rmd x)<+\infty$, where $\|x\|$ denotes the Euclidean norm. We endow $\Pscr_p(\Rbb\p m)$ with the Wasserstein metric $W_p$, that is,

\[
\begin{split}
	W_p(\mu,\nu):= \inf\bigg\{\bigg(\int_{\Rbb\p m\times\Rbb\p m}\|x-y\|\p p&\pi(\rmd x,\rmd y) \bigg)\p{1/p},\\& \pi\in\Pscr_p(\Rbb\p m\times\Rbb\p m):\ \pi(\cdot\times\Rbb\p m)=\mu(\cdot),\  \pi(\Rbb\p m\times\cdot)=\nu(\cdot)\bigg\}.
\end{split}
\] 
From the very definition of $W_p$, we have that for every $\Rbb\p m$-valued  and $p$-integrable random vectors $\et,\xi$ with law $\Lscr(\et)$ and $\Lscr(\xi)$, respectively, the estimate $W_p(\Lscr(\et),\Lscr(\xi))\leq\Ebb[|\et-\xi|\p p]\p{1/p}$ holds.

We consider the following $m$-dimensional mean field SDE:
\begin{equation}\label{eq:mf.eq}
	X_t=\xi+\int_0\p tb(s,X_s,\Lscr(X_s))\rmd s+\int_0\p t\sig(s,X_s,\Lscr(X_s))\rmd B_s,\quad t\in[0,T],
\end{equation}
where $\Lscr(X_s)(=\Pbb_{X_s})$ denotes the law of $X_s$ under $\Pbb$.
\begin{assumption}\label{ass:ass.coeff} We assume that $\Hscr$ is an initial $\sig$-field independent of $\Fscr\p B_\infty$. Moreover, we require
	
	(1) $\xi\in L\p2(\Hscr,\Pbb;\Rbb\p m)$. 
	
	(2) $\|\sig(t,x,\mu)\|_F\p2+\|b(t,x,\mu)\|\p2\leq M(1+\|x\|\p2)$, $(t,x,\mu)\in[0,T]\times\Rbb\p m\times\Pscr_2(\Rbb\p m)$.\newline Here,
	$\|A\|_F=\sqrt{\mathrm{tr}(AA\p\top)}$ denotes the Frobenius norm of the matrix $A\in\Rbb\p{m\times d}$.
 
 (3) Let $h,g,\chi:\Rbb_+\longrightarrow\Rbb_+$ be strictly increasing continuous functions satisfying
  \begin{gather}
  \displaystyle\|\sig(t,x,\mu)-\sig(t,y,\nu)\|_F\leq h(\|x-y\|)+\chi(W_1(\mu,\nu)),\label{eq:ass.coeff.sig}
  		\\
  		\|b(t,x,\mu)-b(t,y,\nu)\|\leq g(\|x-y\|)+\chi(W_1(\mu,\nu)).\label{eq:ass.coeff.b}
   \end{gather} 

 (4) We assume that $h(0+)=\chi(0+)=0$ and $h$ satisfies the Yamada condition condition \[\int_{0}\p1h\p{-2}(r)\rmd r=\infty.\]
 
 (5) There exists a strictly increasing, continuous, concave function $G$ such that $G(0+)=0$ and
 \begin{align}
 	& G(x)\geq g(x)+\frac{m-1}{x}h\p2(x), \quad x>0,\label{eq:co.G1}\\&\int_0\p1\frac1{G(x)}\rmd x=+\infty\label{eq:co.G2}. 
 \end{align} 

\end{assumption}
\begin{remark}\label{rem:onAssumptions}
We observe that, if $m=1$, then the conditions on $G$ are actually conditions on $g$. In this special case, they read: $g$ is a strictly increasing, continuous and concave function such that $g(0+)=0$ and $\int_0\p1\frac1{g(x)}\rmd x=+\infty$. We also remark that, if $b=0$ and $m>1$, then \eqref{eq:co.G2} is equivalent to $\int_0\p1\frac {x}{h\p2(x)}\rmd x=+\infty$, that is stronger than the  Osgood type condition on $h\p2$. A concrete non-Lipschitz example satisfying this assumption is $h(x)= x(\ln\frac{e}{x})\p{1/2}$ for $x\in(0,1)$.
\end{remark}

Let $\mub=(\mu_t)_{t\in[0,T]}\in C([0,T],\Pscr_1(\Rbb\p m))$ be a fixed continuous function of probability measures and let Assumption \ref{ass:ass.coeff} be satisfied. We consider the SDE
\begin{equation}\label{eq:cl.eq}
	X_t\p{\mub}=\xi+\int_0\p tb(s,X_s\p{\mub},\mu_s)\rmd s+\int_0\p t\sig(s,X_s\p{\mub},\mu_s)\rmd B_s,\quad t\in[0,T].
\end{equation}

The following lemma is a consequence of the classical theory of SDEs

\begin{lemma}\label{lem:ex.un}
Let Assumption \ref{ass:ass.coeff} hold. Then, for every $\mub\in C([0,T],\Pscr_1(\Rbb\p m))$, equation \eqref{eq:cl.eq} admits a pathwise unique strong solution $X\p{\mub}$ such that $\sup_{t\leq T} \|X\p{\mub}_t\|\p2\in L\p1(\Pbb)$. Moreover, for every $\mu$, we have $\Ebb[
\|X\p\mu_t\|\p2]\leq C$ and $\Ebb[\|X\p{\mu}_t-X\p\mu_s\|\p2]\leq K|t-s|$, $s,t\in[0,T]$ where $C>0$ and $K>0$ depend only on $T$ and $M$, but neither on $\mu$ nor on $s,t$.
\end{lemma}
\begin{proof}
	By a well known classical result (see \cite{SKO65}, p.\ 59), we have that \eqref{eq:cl.eq} admits a weak solution. Indeed, the coefficients $\widetilde b(t,x):=b(t,x,\mu_t)$ and $\widetilde \sig(t,x):=\sig(t,x,\mu_t)$ are jointly continuous in $(t,x)$. Moreover, by Assumption \ref{ass:ass.coeff}, they satisfy the hypothesis of  \cite[Theorem 2.1]{AS13} (see also \cite[Theorem 4]{YW712}). Hence, \eqref{eq:cl.eq} admits a pathwise unique weak solution.  By the Yamada--Watanabe principle, \eqref{eq:cl.eq} has a strong solution, which is unique in law. The pathwise uniqueness of the strong solution follows by \cite[Theorem 2.1]{AS13}. 
	
	We now come to the estimate of $\sup_{t\leq T} \|X\p{\mub}_t\|\p2\in L\p1(\Pbb)$. This is a classical estimate and we therefore omit details. The idea is to exploit sublinear growth (i.e., Assumption \ref{ass:ass.coeff}(2)) for the coefficients in such a way to apply, for any fixed $R>0$, Gronwall's inequality to the function $u(t)=\Ebb[\sup_{s\leq t\wedge\tau_R} \|X\p{\mub}_s\|\p2]$, where $\tau_R:=\inf\{t\geq0: \|X_t\p{\mub}\|\geq R\}$, and then Fatou's lemma, since $\tau_R\uparrow+\infty$ as $R\rightarrow+\infty$ (as, for example, in the proof of \cite[Theorem 21.11]{Schilling21}). Thanks to  Assumption \ref{ass:ass.coeff}(2), we see that the obtained estimate is not only independent of $t\in[0,T]$ but also of $\mu\in\Pscr_1(\Rbb\p m)$. Once having the uniform estimate of $\|X_t\p{\mub}\|\p2$ in $L\p1$, the $L\p2$-continuity estimate is immediate (see, e.g., the proof of \cite[Corollary 21.30]{Schilling21}). The proof is complete.
\end{proof}
We turn now to the existence theorem. To be more precise, we show that the strong existence of a solution of \eqref{eq:mf.eq} is a consequence of Lemma \ref{lem:ex.un} and of Schauder's fixed-point theorem. 
\begin{theorem}\label{thm:existence}
	Let $\xi$, $b$ and $\sig$ satisfy Assumption \ref{ass:ass.coeff}. Then \eqref{eq:mf.eq} has a strong solution. 
\end{theorem}
Let us begin with some useful remarks that will be used in the proof of this theorem.
\begin{remark}\label{rem:us.prop} 
We discuss some consequences of Assumption \ref{ass:ass.coeff}(4) and introduce some notation that will be directly used in the proof of Theorem \ref{thm:existence} below.
\begin{itemize}
\item[(i)] First we notice that, by Assumption \ref{ass:ass.coeff}(4), there exists a real sequence$(a_n)_n$ such that $a_0:=1$,  $a_n\downarrow0$, $n\rightarrow+\infty$, and $\int_{a_n}\p{a_{n-1}}h\p{-2}(r)\rmd r=n$, $n\geq1$. Hence, for every $n\geq1$, there is a density $\varrho_n$  supported in $(a_n, a_{n-1})$ (i.e., $\int_{a_{n}}\p{a_{n-1}}\rho_n(r)\rmd r=1$) such that $0\leq\rho_n(r)\leq\frac2{nh\p2(r)}$. We define the function $\psi_n\in C\p2 (\Rbb)$
\[
\psi_n(x)=\int_0\p{x}\int_0\p y\rho_n(r)\rmd r\,\rmd y,\quad x\geq0.
\]
Obviously,  $\psi_n$ is an increasing function such that $\max\{x-a_{n-1},0\}\leq\psi_n(x)\leq x$ and $\psi_n(x)\uparrow x$,  $n\rightarrow+\infty$, $x\geq0$. Moreover, $\psi_n(x)=0$ for $x\in[0,a_n]$, $|\psi\p\prime_n(x)|\leq1$ and $\psi_n\p{\prime\prime}(x)=\rho_n(x)$, $x\geq0$.

\item[(ii)] We then define the function $f_n(x):=\psi_n(\|x\|)$, $n\geq1$, $x\in\Rbb\p m$. Since $x\mapsto\|x\|$ is convex and increasing, we deduce that $f_n=\psi_n\circ\|\cdot\|$ is convex too, $\psi_n$ being convex. Clearly, $f_n\in C\p2(\Rbb\p m)$ and $f_n(x)=0$ if $\|x\|\leq a_n$.
By direct computation it is easy to see that $\nabla f_n(x)=\psi_n\p\prime(\|x\|)\frac x{\|x\|}$ and
\[
\partial\p2_{x_ix_j}f_n(x)=\psi_n\p{\prime\prime}(\|x\|)\frac{x_ix_j}{\|x\|\p2}-\psi_n\p{\prime}(\|x\|)\bigg(\frac{x_ix_j}{\|x\|\p3}-1_{\{i=j\}}\frac1{\|x\|}\bigg)
\]
 Denoting by $H_n(x)$ the Hessian matrix of $f_n$ at $x\in\Rbb\p m$, taking the trace, denoted by $\tr$, and using that $\psi_n\p\prime(\|x\|)\leq1_{\{\|x\|\geq a_n\}}$,  yields
\begin{equation}\label{eq:tr.est}
\begin{split}
\tr\{H_n(x)\}=\rho_n(\|x\|)+\frac{m-1}{\|x\|}\psi_n\p\prime(\|x\|)
&\leq
\rho_n(\|x\|)+\frac{m-1}{\|x\|}1_{\{\|x\|\geq a_n\}}
\\&
\leq\rho_n(\|x\|)+\frac{m-1}{\|x\|}1_{\{\|x\|\neq 0\}}.
\end{split}
\end{equation}
\item[(iii)] Since $f_n$ is convex, we have that $H_n(x)$ is positive semi-definite. So, from \cite[Theorem 1]{Coope93} (see also \cite[Remark 2.20 (c)]{AS13}), for any other $A\in\Rbb\p{m\times d}$, we have 
\[
\tr\{H_n(x)AA\p\top\}\leq\tr\{H_n(x)\}\tr\{AA\p\top\}=\tr\{H_n(x)\}\|A\|_F\p2.
\]
\end{itemize}
\end{remark}

\begin{proof}[Proof of Theorem \ref{thm:existence}]
We consider the space $V:=C([0,T],\Pscr_1(\Rbb\p m))$ endowed with the metric 
\[
d(\mu,\nu):=\sup_{t\in[0,T]}W_1(\mu_t,\nu_t),\quad \mu,nu\in V,
\]
and observe that $(V,d)$ is a Banach space.
By the moment estimates in Lemma \ref{lem:ex.un}, we know that $\Lscr(X\p\mu_t)\in\Pscr_2(\Rbb\p m)\subseteq\Pscr_1(\Rbb\p m)$. Moreover, the continuity estimate in Lemma \ref{lem:ex.un} yields that $(\Lscr(X\p\mu_t))_{t\in[0,T]}\in C([0,T],\Pscr_1(\Rbb\p m))$. Indeed,
\[
W_1(\Lscr(X\p\mu_t),\Lscr(X\p\mu_s))\leq W_2(\Lscr(X\p\mu_t),\Lscr(X\p\mu_s))\leq\Big(\Ebb\big[\|X_t\p\mu-X\p\mu_s\|\p2\big]\Big)\p{1/2}\leq \sqrt K|t-s|\p{1/2},\quad  s,t\in[0,T].
\]
Hence,  we can define a mapping $\Phi: (V,d)\longrightarrow (V,d)$  by putting $\Phi(\mu)_t=\Lscr(X_t\p\mu)$, $t\in[0,T]$. We show that $\Phi$ is continuous (in the topology induced by the metric $d$). To this aim, we  adapt the original argument by Yamada and Watanabe in \cite[Theorem 1]{YW711}, as in the proof of \cite[Theorem 2.1]{AS13}.

Let $\mub\p{i}\in V$ and consider the associated pathwise unique strong solution $X\p{\mub\p i}$ to \eqref{eq:cl.eq},  $i=1,2$. We define
\begin{equation}\label{eq:D}
D_t:=D_t(\mub\p1,\mub\p2):=X_t\p{\mub\p 1}-X_t\p{\mub\p 2},\quad t\in[0,T].
\end{equation}
We also introduce the notation
\[
\overline\sigma_t:=\sig(t,X\p{\mub\p 1}_t,\mu\p 1_t)-\sig(t,X\p{\mub\p 2}_t,\mu\p 2_t),\qquad\overline b_t:=b(t,X\p{\mub\p 1}_t,\mu\p 1_t)-b(t,X\p{\mub\p 2}_t,\mu\p 2_t),\quad t\in[0,T].
\]
Also observe that, by Assumption \ref{ass:ass.coeff}(2) and Lemma \ref{lem:ex.un}, the mappings 
\[(\om,t)\mapsto \|b(t,X\p{\mub\p i}_t(\om),\mu\p i_t)\|\qquad \textnormal{and}\qquad  (\om,t)\mapsto \|\sig(t,X\p{\mub\p i}_t(\om),\mu\p i_t)\|\p2_F\] are $\Pbb(\rmd\om)\rmd t$-integrable. We now make use of the facts summarized in Remark \ref{rem:us.prop}, without further mention. First we apply the It\^o formula to the semimartingale $\psi_n(\|D\|)=f_n(D)$ and observe that, since $\|\nabla f_n\|\leq1$, by the Cauchy-Schwartz inequality $\|x\p\top A\|\leq\|x\|\|A\|_F$ for $x\in\Rbb\p m$ and $A\in\Rbb\p{m\times d}$, the stochastic integrals with respect to the components of $W$ are true martingales. Taking now the expected value and applying It\^o's formula and Fubini's theorem, Assumption \ref{ass:ass.coeff} yields

\begin{equation}\label{eq:ito.psi.d}
	\begin{split}
\Ebb[\|D_t\|]-a_{n-1}\leq\Ebb\big[f_n(D_t)\big]&=\Bigg|\int_0\p t\Ebb\bigg[\nabla f_n(D_s)\ol b_s+ \frac12tr\{ H_n(D_s)\ol\sig_s\ol\sig_s\p\top\}\bigg]\rmd s\Bigg|
\\&\leq
\int_0\p t\Ebb\bigg[\|\ol b_s\|+\frac12\tr\{H_n(D_s)\}\|\ol\sig_s\|_F\p2\bigg]\rmd s
\\&\leq
\int_0\p t\Ebb\bigg[\|\ol b_s\|+\frac12\bigg(\rho_n(\|D_s\|)+\frac{m-1}{\|D_s\|}1_{\{\|D_s\|\geq a_n\}}\bigg)\|\ol\sig_s\|_F\p2\bigg]\rmd s.
\end{split}
\end{equation}
We now estimate the second summand in the expected value on the right-hand side of \eqref{eq:ito.psi.d} separately.
By Remark \ref{rem:us.prop} and \eqref{eq:ass.coeff.sig}, we have
\begin{equation}\label{eq:est.sigF}
\begin{split}
\bigg(\rho_n&(\|D_s\|)+\frac{m-1}{\|D_s\|}{1}_{\{\|D_s\|\geq a_n\}}\bigg)\|\bar{\sigma_s}\|_F^2
\\&\leq
2\bigg(\rho_n(\|D_s\|)+\frac{m-1}{\|D_s\|}{1}_{\{\|D_s\|\geq a_n\}}\bigg)\big(h^2(\|D_s\|)+\chi^2\big(d(\mu\p1,\mu\p2)\big)\big)
\\&\leq
2\bigg(\frac{2}{nh\p2(\|D_s\|)}1_{(a_n,a_{n-1})}(\|D_s\|)+\frac{m-1}{\|D_s\|}{1}_{\{\|D_s\|\geq a_n\}}\bigg)\big(h^2(\|D_s\|)+\chi^2\big(d(\mu\p1,\mu\p2)\big)\big)
\\&\leq
2\bigg(\frac2n+(m-1)\frac{h\p2(\|D_s\|)}{\|D_s\|}1_{\{\|D_s\|\neq0\}}\bigg)+2\bigg(\frac{2}{nh\p2(a_n)}+\frac{m-1}{a_n}\bigg)\chi^2\big(d(\mu\p1,\mu\p2)\big)
\end{split}
\end{equation}
Plugging the estimate \eqref{eq:est.sigF} in \eqref{eq:ito.psi.d} and using \eqref{eq:ass.coeff.b} and \eqref{eq:co.G1}, now yields
\begin{equation}\label{eq:ito.psi}
	\begin{split}
\Ebb[\|D_t&\|]-a_{n-1}
\leq
T\bigg(\frac2n+\chi(d(\mu_1,\mu_2))+\Big(\frac{2}{nh^2(a_n)}+\frac{m-1}{a_n}\Big)\chi^2(d(\mu\p1,\mu\p2))\bigg)
\\
&\quad+\int_0^T \Ebb\bigg[g(\|D_s\|)+(m-1)\frac{h^2(\|D_s\|)}{\|D_s\|}{1}_{\{\|D_s\|\neq0\}}\bigg]\rmd s
\\
&\leq
T\bigg(\frac2n+\chi(d(\mu\p1,\mu\p2))+\Big(\frac{2}{nh^2(a_n)}+\frac{m-1}{a_n}\Big)\chi^2(d(\mu\p1,\mu\p2))\bigg)
\\
&\hspace{8cm}+\int_0^t G(\Ebb[\|D_s\|])\rmd s,
\end{split}
\end{equation}
where in the last step we have used Jensen's inequality combined with the concavity of $G$ (Assumption \ref{ass:ass.coeff}(5)). Then, by taking the limit $d(\mub\p1,\mub\p2)\rightarrow0$ and putting $u(t):=\limsup_{d(\mub\p1,\mub\p2)\rightarrow0}\Ebb\big[\|D_t(\mu\p1,\mu\p2)\|\big]$, we obtain
\[
	u(t)\leq T\frac2n+\int_0\p TG\big(u(s))\rmd s,\qquad t\in[0,T].
\]
Consequently, letting tend $n\rightarrow+\infty$ implies 
\begin{equation}\label{eq:est.u}
	u(t)\leq \int_0\p TG\big(u(s))\rmd s,\qquad t\in[0,T].
\end{equation}
Note that $u(0)=0$. To apply the  Bihari--La Salle inequality (see \cite[Lemma 3.6]{Mao95}), we must show that $t\mapsto u(t)$ is continuous. We observe that, by Lemma \ref{lem:ex.un}, the estimate 

\begin{equation}\label{eq:equic}
	\Ebb[\|X\p{\mu}_t-X\p\mu_s\|\p2]\leq K|t-s|,\quad t,s\in[0,T],
	\end{equation}
 holds, where $K>0$ is a constant depending on $T$ and on $M$, but not on $\mu$ nor on $s,t\in[0,T]$. So we have 
\begin{equation}\label{eq:unif.cont}
	\begin{split}
\Ebb[\|D_t-&D_s\|]
\leq\Big(\Ebb\big[\|D_t-D_s\|\p2\big]\Big)\p{1/2}
\\&\leq\sqrt 2\bigg(\Big(\Ebb\Big[\|X\p{\mu\p1}_t-X\p{\mu\p1}_s\|\p2\Big]\Big)\p{1/2}+\Big(\Ebb\Big[\|X\p{\mu\p2}_t-X\p{\mu\p2}_s\|\p2\Big]\Big)\p{1/2}\bigg)\leq \,2\sqrt {2K}|t-s|\p{1/2}.
\end{split}
\end{equation}
Hence,
\[
	|u(t)-u(s)|\leq\limsup_{d(\mu\p1,\mu\p2)\rightarrow0}\Ebb[\|D_t(\mu\p1,\mu\p2)-D_s(\mu\p1,\mu\p2)\|]\leq 2\sqrt{2K}|t-s|\p{1/2}\longrightarrow0, \quad\textnormal{ as }\quad s\rightarrow t.
\] 
Therefore, we can apply the Bihari--La Salle inequality, and by Assumption \ref{ass:ass.coeff}(5), we get $u(t)=0$, $t\geq0$.
But this now shows the continuity of $\Phi$. Indeed, by \eqref{eq:ito.psi} we get
\begin{equation}\label{eq:cont.Phi}
	\begin{split}
	0\leq\limsup_{d(\mu\p1,\mu\p2)\rightarrow0}d(\Phi(\mu\p1),\Phi(\mu\p2))
	&\leq
	\limsup_{d(\mu\p1,\mu\p2)\rightarrow0}\sup_{t\leq T}\Ebb[\|D_t(\mu\p1,\mu\p2)\|]
	\\&
	\leq\Big(a_{n-1}+\frac{2T}n\Big)+\int_0\p TG(u(s))\rmd s
	\\&=
	\Big(a_{n-1}+\frac{2T}n\Big)\longrightarrow0,\quad \textnormal{ as }\ n\rightarrow+\infty.
\end{split}	
\end{equation}
The proof of the continuity of $\Phi:(V,d)\longrightarrow (V,d)$ is complete.

Let us now introduce the space $\Kscr:=\{\mu\in\Pscr_1(\Rbb\p m): \int_\Rbb\|x\|\p2\mu(\rmd x)\leq C\}$, where $C>0$ is the same as in Lemma \ref{lem:ex.un}.  We show that $\Kscr$ is a compact subset of $(\Pscr_1(\Rbb\p m),W_1)$. For this end, consider any $(\mu_n)_n\subseteq\Kscr$ and let $\xi_n$ be a random vector on $(\Om,\Fscr,\Pbb)$ obeying the law $\mu_n$, for all $n\geq1$. Then, by the Chebyshev--Markov inequality, we get with choosing $C_\ep:=C\p{1/2}/\ep$, $\ep>0$,
\[
\mu_n(\Rbb\p m\setminus\{\|x\|\leq C_\varepsilon\})=\Pbb[\|\xi_n\|> C_\ep]\leq\frac{\Ebb[\|\xi_n\|]}{C_\ep}\leq\frac{\Ebb[\|\xi_n\|\p2]\p{1/2}}{C_\ep}\leq\ep,\quad n\geq1,
\] 
meaning that $(\mu_n)_n$ is tight. Thus, due to the Prohorov Theorem, there exists a subsequence $(\mu_{n_k})_k$ and $\mu\in\Pscr(\Rbb\p m)$ such that $\mu_{n_k}\rightharpoonup\mu$ weakly, as $k\rightarrow+\infty$. By the definition of $\Kscr$, $\sup_n\int_{\Rbb\p m}\|x\|\p2\mu_n(\rmd x)\leq C$. Thus, the sequence $(\mu_n)_n$ is uniformly integrable:
\[
\int_{\{\|x\|>R\}}\|x\|\mu_n(\rmd x)\leq \int_{\{\|x\|>R\}}\|x\|\frac{\|x\|}{R}\mu_n(\rmd x)\leq \frac{C}{R}\longrightarrow0, \quad\textnormal{ as }\ R\rightarrow +\infty.
\] 
Consequently, by \cite[Theorem 7.2]{V21}, we get $W_1(\mu_{n_k},\mu)\longrightarrow0$, as $k\rightarrow+\infty$.
It remains to show that $\mu\in\Kscr$. But this is obvious as, by the above weak convergence, for every $N\geq1$, we have 
\[
\int_{\Rbb\p m}((\|x\|\p2\wedge N)\mu(\rmd x)=\lim_{n\rightarrow+\infty}\int_{\Rbb\p m}(\|x\|\p2\wedge N)\mu_n(\rmd x)\leq C.
\] 
Now compactness follows by monotone convergence. 

By Lemma \ref{lem:ex.un} and \eqref{eq:equic}, we have $\Phi:\Escr\longrightarrow\Escr$, where $\Escr$ is defined by
\[
\Escr:=\{\mu\in C([0,T],\Kscr): \mu_0=\Pbb_\xi\ \textnormal{ and } W_1(\mu_t,\mu_s)\leq \sqrt K|t-s|\p{1/2}\},
\]
with the same constant  $K>0$ appearing in \eqref{eq:equic}.
Obviously, $\Escr$ is a closed, convex subspace of $\Kscr$-valued equi-continuous functions. Thus, by Ascoli--Arzel\`a Theorem, since $\Kscr$ is compact, $\Escr$ is relatively compact in $(V,d)$. Hence, its closure $\overline\Escr$ is a compact subset of $(V,d)$. But as $\Escr$ is a closed subset of $(V,d)$, it is compact. 

We now consider the space $\Mscr_1(\Rbb\p m)$ of finite signed measures such that $\int_\Rbb\|x\||\mu|(\rmd x)<+\infty$, where  $|\mu|$ denotes the total variation of a finite-valued signed measure $\mu$. Then $\Mscr_1(\Rbb\p m)$ endowed with the Kantorovic--Rubinstein norm,
\[
\|\mu\|_{\mathrm{KR}}:=\sup_{f\in\mathrm{Lip_{1,0}}}\int_{\Rbb\p m} f(x)\mu(\rmd x)
\]
is a normed vector space, where $\mathrm{Lip_{1,0}}$ denotes the class of Lipschitz functions on $\Rbb\p m$ with Lipschitz constant $L=1$ and vanishing in $x=0$. Obviously,  the set $W:=C([0,T],\Mscr_1(\Rbb\p m))$  endowed with the norm $\|\mu\|_{\mathrm{KR},\infty}:=\sup_{t\leq T}\|\mu_t\|_{\mathrm{KR}}$ is a normed vector space, and so it is a Hausdorff locally convex topological vector space. Clearly, $V\subseteq W$ and, by the Kantorovic--Rubinstein duality, 
\[
d(\mu,\nu)=\sup_{t\leq T} W_1(\mu_t,\nu_t)=\sup_{t\leq T}\|\mu_t-\nu_t\|_{\mathrm{KR}}=\|\mu-\nu\|_{\mathrm{KR},\infty}, \quad \mu,\nu\in V.
\]
Thus, $(V,d)$ is continuously embedded in $(W,\|\cdot\|_{\mathrm{KR},\infty})$.

Summarizing, we have shown that $\Phi$ is continuous on $(V,d)\subseteq (W,\|\cdot\|_{\mathrm{KR},\infty})$ and maps the compact subset $(\Escr,d)\subseteq(V,d)\subseteq (W,\|\cdot\|_{\mathrm{KR},\infty})$ into itself, i.e., $\Phi(\Escr)\subseteq\Escr$. From Schauder's fixed-point theorem, it follows that $\Phi$ has a fix-point in $\Escr$: There exists $\mu\in\Escr$ such that $\mu_t=\Phi(\mu)_t=\Lscr(X\p\mu_t)$, $t\in[0,T]$. This shows the existence of a strong solution of \eqref{eq:mf.eq}. The proof of the theorem is complete.
\end{proof}

\section{Uniqueness of the Strong Solution}\label{sec:uniqueness}
In order to show the strong uniqueness of the solution we strengthen the assumptions on the coefficients:
\begin{theorem}\label{thm:un} Let Assumption \ref{ass:ass.coeff} hold. Suppose further that there exists an increasing, continuous and concave function $q$ such that $q(0)=0$, $\int_{0+}q\p{-1}(r)\rmd r=+\infty$ and

		\begin{align}
	&\displaystyle\|\sig(t,x,\mu)-\sig(t,y,\nu)\|_F\p2\leq q\big(\|x-y\|\p2+W_2\p2(\mu,\nu)\big),\quad t\in[0,T],\  x,y\in\Rbb\p m,\ \mu,\nu\in \Pscr_2(\Rbb\p m),\label{eq:ass.coeff.un.sig}
			\\
		&\Ebb\big[(b(t,\et,\Pbb_\et)-b(t,\zt,\Pbb_\zt))(\et-\zt)\big]\leq q\big(2\Ebb\big[(\et-\zt)\p 2\big]\big), \quad t\in[0,T],\ \et,\zt\in L\p2(\Om,\Fscr_t,\Pbb).	\label{eq:ass.coeff.un.b}
		\end{align}
	Then the strong solution of \eqref{eq:mf.eq} is pathwise unique.
\end{theorem}
\begin{proof} Let $X\p i$, $i=1,2$, be two strong solutions of \eqref{eq:mf.eq} and $D_t=D(\Lscr(X\p1_t),\Lscr(X\p2_t))$ be defined as in \eqref{eq:D}. Then,
\[
\begin{split}
\Ebb[\|D_t\|\p2]&
=
\int_0\p t\Ebb\Big[\big\|\sig(s,X\p1_s,\Lscr(X_s\p1))-\sig(s,X\p2_s,\Lscr(X_s\p2))\big\|\p2\Big]\rmd s
\\&\qquad+2\int_0\p t\Ebb\Big[\Big\langle D_s,b(s,X\p1_s,\Lscr(X_s\p1))-b(s,X\p2_s,\Lscr(X_s\p2))\Big\rangle\Big]\rmd s
\\&
\leq 3 \int_0\p t q\big(2\Ebb\big[\|D_s\|\p2\big]\big)\rmd s, \qquad t\in[0,T].
\end{split}
\]
Since by Lemma \ref{lem:ex.un} the function $t\mapsto u(t)=\Ebb[\|D_t\|\p2]$ is continuous and $u(0)=0$, the Bihari--La Salle inequality yields that $u(t)=0$, for every $t\in[0,T]$. Therefore, $X\p 1$ and $X\p 2$ are indistinguishable, since they are continuous. The proof is complete.
\end{proof}
\begin{remark}\label{rem:coef}
	(i) We notice that, if \eqref{eq:ass.coeff.un.sig} holds with $W_1$ instead of $W_2$, then \eqref{eq:ass.coeff.un.sig}  trivially holds. Moreover, \eqref{eq:ass.coeff.sig} with $h(x)=\chi(x)=q(x^2)$ is verified.
	
	(ii) Observe that \eqref{eq:ass.coeff.un.b} weakens the monotonicity condition
	\[
	\Ebb\big[\langle b(t,\et,\Pbb_\et)-b(t,\zt,\Pbb_\zt),\et-\zt\rangle\big]\leq C\Ebb\big[\|\et-\zt\|\p 2\big], \quad t\in[0,T],\ \ \et,\zt\in L\p2(\Om,\Fscr_t,\Pbb),
	\]
	introduced and studied in \cite[Section 3]{HLL25}. Moreover, \eqref{eq:ass.coeff.un.b} trivially holds if
	\[
	\langle b(t,x,\mu)-b(t,y,\nu),x-y\rangle\leq q\big(\|x-y\|\p2+W_2\p2(\mu,\nu)\big),\quad t\in[0,T],\quad x,y\in\Rbb\p m,\quad \mu,\nu\in \Pscr_2(\Rbb\p m).
	\]
	
	(iii) If $b(t,x,\Pbb_\et)=b_1(t,x)+b_2(t,\Ebb[\et])$ for $t\in[0,T]$, $x\in\Rbb\p m$, $\et\in L\p2(\Om,\Fscr,\Pbb)$, then \eqref{eq:ass.coeff.un.b} reduces to the classical monotonicity condition
	\[
	\langle b_i(t,x_i)-b_i(t,y_i),x_i-y_i\rangle\leq q(|x_i-y_i|\p2),\quad t\in[0,T],\ x_i,y_i\in\Rbb\p m,\ i=1,2.
	\]
	\end{remark}
\bibliographystyle{plainnat}
\bibliography{bibliography}
\end{document}